\documentclass[12pt,leqno]{amsart}
\usepackage[margin=1.5in]{geometry}
\usepackage[T1]{fontenc}
\usepackage{lmodern}
\usepackage{microtype}
\usepackage{amsmath,amssymb,mathtools,esint}
\usepackage{mathrsfs}
\usepackage{amsfonts}
\usepackage[dvipsnames]{xcolor}
\usepackage[hidelinks]{hyperref}
\hypersetup{
 pdftitle={A singular minimizer of a uniformly convex Hessian integral},
 pdfauthor={Arunima Bhattacharya and Anna Skorobogatova},
 pdfsubject={Scalar Hessian integrals and singular minimizers}
}
\newtheorem{theorem}{Theorem}[section]
\newtheorem{proposition}[theorem]{Proposition}
\newtheorem{lemma}[theorem]{Lemma}
\newtheorem{corollary}[theorem]{Corollary}
\theoremstyle{definition}

\theoremstyle{remark}
\newtheorem{remark}[theorem]{Remark}
\numberwithin{equation}{section}
\newcommand{\R}{\mathbb R}
\newcommand{\Sym}{\operatorname{Sym}}
\newcommand{\tr}{\operatorname{tr}}
\newcommand{\supp}{\operatorname{supp}}
\newcommand{\sgn}{\operatorname{sgn}}

\newcommand{\Sing}{\operatorname{Sing}}
\newcommand{\loc}{\mathrm{loc}}

\newcommand{\cT}{\mathcal T}

\newcommand{\ip}[2]{\langle #1,#2\rangle}

\newcommand{\op}{\mathrm{op}}
\newcommand{\dd}{\,\mathrm{d}}
\newcommand{\Id}{\mathrm{Id}}
\newcommand{\divdiv}{\operatorname{div}\operatorname{div}}
\usepackage{comment}

\newcommand{\eps}{\varepsilon}

\newcommand{\ttimes}{%
  \mathbin{\vcenter{\hbox{%
    \ooalign{%
      \hfil\raise.18ex\hbox{$\times$}\hfil\cr
      \hfil\lower.18ex\hbox{$\times$}\hfil\cr
    }%
  }}}%
}

\newcommand{\Tcal}{\mathcal{T}}

\newcommand{\Sscr}{\mathscr{S}}

\newcommand{\Sbb}{\mathbb{S}}

\allowdisplaybreaks[2]

\title[A singular minimizer of a Hessian-dependent functional]{A singular minimizer of a uniformly convex Hessian-dependent functional}
\author{Arunima Bhattacharya}
\address{Department of Mathematics, University of North Carolina at Chapel Hill,
 Chapel Hill, NC 27599, USA}
\email{arunimab@unc.edu}
\author{Anna Skorobogatova}
\address{Institute for Theoretical Studies, ETH Z\"urich, Z\"urich, Switzerland}
\email{anna.skorobogatova@eth-its.ethz.ch}
\date{}

\begin{document}
\begin{abstract}
We construct a smooth, uniformly convex Hessian-dependent energy functional, which admits a singular $C^{1,1}$ minimizer on $\mathbb{R}^5$ that is a two-homogeneous normalization of Cartan's five-dimensional cubic. 
Our construction relies on constructing a quadratic null Lagrangian for the Hessian that is suitably non-degenerate on the image of the Hessian for the candidate minimizer.
This complements a partial regularity result previously obtained by the authors in this setting, and demonstrates that one cannot hope for a fourth-order version of classical second-order elliptic regularity theory \`{a} la De Giorgi-Nash-Moser.

\end{abstract}
\maketitle

\section{Introduction}

In this article, we consider smooth uniformly convex functionals
of the form
\begin{equation}\label{eq:energy}
 \mathcal E_F(v;\Omega)=\int_\Omega F(D^2v)\dd x,
 \qquad v\colon\Omega\subset\R^n\longrightarrow\R
\end{equation}
on a bounded domain $\Omega$, with $F\in C^\infty(\Sym(n))$ satisfying
\begin{equation}\label{eq:ellipticity}
 \lambda|N|^2\leq D^2F(M)[N,N]\leq\Lambda|N|^2
 \qquad\text{for all }M,N\in\Sym(n),
\end{equation}
for some $0<\lambda\leq\Lambda<\infty$. The associated Euler--Lagrange equation satisfied by critical points is
$\divdiv DF(D^2v)=0$, understood weakly as
\begin{equation}\label{eq:EL}
 \int_\Omega DF(D^2v) : D^2\varphi\dd x=0
 \qquad\text{for every }\varphi\in C_c^\infty(\Omega).
\end{equation}
The upper bound in \eqref{eq:ellipticity} gives quadratic growth of $F$
and linear growth of $DF$, so both the functional and its first variation
are well-defined on $W^{2,2}(\Omega)$. Moreover, uniform convexity of $\mathcal{E}_F$ implies that a
weak solution of \eqref{eq:EL} is the unique minimizer in
$v+W^{2,2}_0(\Omega)$.

Hessian-dependent energies and their associated fourth-order Euler--Lagrange equations occur naturally in elasticity and models of folding
\cite{Tou62,Mindlin64,JK00,FJM06}, in liquid crystals \cite{AG87}, in the affine Plateau problem \cite{TW05}, and in the study of canonical metrics that are minimizers or critical points of various curvature-dependent energies \cite{Wil65,Don02,P,GV15}. Our particular motivation comes from Hamiltonian stationary Lagrangian
submanifolds, see for instance \cite{Oh90,Oh93} for the variational formulation, and \cite{SW01} for the constrained area
minimization problem. Recall that the volume of a Lagrangian gradient graph $\{(x,Dv(x)): x\in \Omega\}\subseteq \mathbb{R}^n\times \mathbb {R}^n$ is
\begin{equation}\label{e:HamStat}
 \int_\Omega\sqrt{\det\bigl(I+(D^2v)^2\bigr)}\dd x,
\end{equation}
and the vanishing of its first variation under compactly supported Hamiltonian variations gives the Hamiltonian stationary equation for the potential $v$ (see \cite[eqn (1.1)]{CW19}). Note that the energy \eqref{e:HamStat} is smooth and uniformly convex when $|D^2 v|$ is sufficiently small, but is not globally uniformly convex on $\Sym(n)$, see e.g. \cite{BOS26} for a more detailed discussion.

For gradient-dependent energy functionals of the form $\int_\Omega\Phi(Dv)\dd x$, smoothness of scalar minimizers under the corresponding assumptions follows from 
De Giorgi--Nash theory \cite{DG57,Nash58}, applied to the second order Euler--Lagrange
equation. Equation \eqref{eq:EL} is instead fourth order, and this
argument is no longer available. This raises the natural question of whether one expects scalar minimizers of Hessian-dependent energies \eqref{eq:energy} to also be smooth. Although a Schauder theory was developed by Bhattacharya \& Warren in \cite{BW21}, demonstrating that $C^{2,\alpha}$ critical points of \eqref{eq:energy} are smooth, this left open the question of whether $W^{2,2}$ critical points are $C^{2,\alpha}$. When $n=2$, this is indeed true, as demonstrated by Ne\v{c}as \cite[Théorème 3]{Nec67}.

In general dimensions, merely $\eps$-regularity results for $W^{2,\infty}$ critical points were known; see \cite{Bha22} for the case of small $L^\infty$ Hessian excess (or oscillation), and our earlier work \cite{BS25} for the case of small $L^p$ Hessian excess, i.e., small BMO seminorm. The latter in particular yields a partial regularity result for bounded-Hessian critical points: there exists $p_0\in(2,3)$ such that the critical point $v$ is smooth outside the relatively closed set
\begin{equation}\label{e:singset}
 \Sigma_v=\left\{x\in\Omega:
 \liminf_{r\downarrow0}\fint_{B_r(x)}
 |D^2v-(D^2v)_{B_r(x)}|^{p_0}\dd y>0\right\},
\end{equation}
which satisfies the dimension bound
\begin{equation}\label{e:dim}
 \dim_{\mathcal H}\Sigma_v\leq n-p_0.
\end{equation}
Here $(h)_B=\fint_Bh$ denotes the average over $B$.

In this article, we show that singular minimizers for energies of the form \eqref{eq:energy} do exist in general. More precisely, we construct a uniformly convex energy functional of the form \eqref{eq:energy} on the space $\Sym(5)$ of symmetric $5\times 5$ matrices, for which the associated (unique) scalar minimizer on $\R^5$ is a 2-homogeneous function that is not a quadratic polynomial, and thus even lies in $W^{2,\infty}$, but not in $C^2$ at the origin. This provides a first example of a singular minimizer in this setting.

Before providing the precise statement of our main theorem, we recall the Cartan isoparametric cubic $P$ on $\R^5$ (see e.g. \cite{Car39}), which is the three-homogeneous polynomial
\begin{equation}\label{e:Cartan-cubic}
     P(x)=x_1^3+\frac{3x_1}{2}
  \bigl(x_3^2+x_4^2-2x_5^2-2x_2^2\bigr)
  +\frac{3\sqrt3}{2}
  \bigl(x_2x_3^2-x_2x_4^2+2x_3x_4x_5\bigr),
\end{equation}
where $x= (x_1,\dots, x_5)$. We then consider the two-homogeneous extension of its link in $\mathbb{S}^4$, i.e.
\begin{equation}\label{eq:potential}
u(x)=\begin{cases}\frac{P(x)}{|x|},&x\neq0,\\0,&x=0.\end{cases}
\end{equation}
Note that $u$ is smooth away from the origin and has uniformly bounded Hessian on $\R^5 \setminus \{0\}$ thanks to the zero-homogeneity of the latter, but $u$ is not twice differentiable at the origin. Indeed, if it were, its two-homogeneity would force it to agree with its quadratic Taylor polynomial, but $u$ is odd. In particular, it is easy to check that $0 \in \Sigma_u$ as defined in \eqref{e:singset}.

Our main result is the following, in which we demonstrate the existence of an energy functional of the form \eqref{eq:energy} for which the above function $u$ is a global minimizer. 
\begin{theorem}\label{thm:main}
There exists $F\in C^\infty(\Sym(5))$ satisfying
\eqref{eq:ellipticity} for a suitable choice of ellipticity constants $\lambda, \Lambda$, such that the function $u$ in
\eqref{eq:potential} is a weak solution of \eqref{eq:EL} on $B_1\subset \R^5$. Moreover, for every $\varphi\in W^{2,2}_0(B_1)$,
\begin{equation}\label{eq:strictminimum}
 \mathcal E_F(u+\varphi;B_1)-\mathcal E_F(u;B_1)
 \geq\frac\lambda2\int_{B_1}|D^2\varphi|^2\dd x.
\end{equation}
In particular, $u$ is the unique minimizer with the smooth data
$u=P$ and $\partial_\nu u=2P$ on $\partial B_1$.
\end{theorem}

For $u\in W^{2,\infty}(B_1)$, we will henceforth let
\begin{align*}
    &\mathrm{Reg}(u) := \{x: \text{there exists a neighborhood $U \ni x$ with $u \in C^2(U)$}\}, \\
    &\mathrm{Sing}(u) := B_1 \setminus \mathrm{Reg}(u).
\end{align*}
Note that $\Sing(u)$ is closed and is contained in $\Sigma_u$ as defined in \eqref{e:singset}. As a consequence of Theorem \ref{thm:main}, we easily produce examples with codimension 5 singular sets in any dimension.
\begin{corollary}\label{cor:dimensions}
For each $n\geq5$, there exists a smooth integrand $F\colon\Sym(n)\to\R$
satisfying \eqref{eq:ellipticity} for some $\lambda(n), \Lambda(n)$, and
a minimizer $u\in W^{2,\infty}(B_1)$ such that
\begin{equation}\label{e:higher-dim-spine}
 \Sing(u) = \Sigma_u =(\{0\}\times\R^{n-5})\cap B_1.
\end{equation}
In particular, $\dim_{\mathcal H}\Sing(u)=n-5$.
\end{corollary}

\begin{remark}\label{r:strict-conv}
    Note that strict convexity alone is insufficient for an $\eps$-regularity theorem. This stems from the fact that, unlike in the uniformly convex case, there is no longer a Schauder theory for $C^{2,\alpha}$-regular critical points in general. Indeed, for each even integer $p \geq 4$, an elementary one-dimensional example is
    \[
        F_p(m)= m^p, \qquad v(t)=\sgn(t)|t|^{2+\frac{1}{p-1}},
    \]
    which corresponds to the $p$-biLaplacian; cf. \cite[Remark 1.1]{BB18} (see also \cite{Moo20}) for a corresponding class of examples for the $p$-Laplacian. It is easy to check that $F$ is smooth and strictly convex but not uniformly convex, and $v$ lies in the class $C^{2,\frac{1}{p-1}}$ but not in $C^3$.
\end{remark}

When combined with \eqref{e:dim}, this leaves open the question of the optimal $p\in [p_0,5]$ such that $\dim_\mathcal{H} \Sigma_u \leq n-p$ for any minimizer $u$ of an energy functional of the form \eqref{eq:energy}.

One may compare our partial regularity results to those for vectorial minimizers of gradient-dependent energies. Indeed, one may view \eqref{eq:energy} as a gradient-dependent energy on gradient vector fields. The regularity theory for vectorial critical points $u$ of uniformly elliptic gradient-dependent integrands has been the subject of extensive study; the same dimension estimate \eqref{e:dim} on the singular set is known for some $p_0\in (2,3)$ in that setting \cite{GiMi68, Giusti1969Singolarita, Morrey66, GG78}; while singular minimizers of uniformly convex
first-gradient integrals do indeed exist. The first result in the latter direction is due to Ne\v{c}as \cite{Nec77}, where a
smooth, uniformly convex gradient-dependent energy functional admitting a singular minimizer is constructed in sufficiently high dimension. In subsequent work \cite{HLN96}, Hao, Leonardi \& Ne\v{c}as produced an example with an isolated singularity in dimension five. \v{S}ver\'ak \& Yan \cite{SY00} then proceeded to demonstrate the existence of a
Lipschitz (see also \cite{SY02} for non-Lipschitz examples) singular minimizer of a suitable gradient-dependent energy from $\R^3$ to $\R^5$ via first constructing a suitable null Lagrangian, and then finding a uniformly convex function whose derivative matches that of the null Lagrangian on gradients. More recently, Mooney \& Savin \cite{MS16} obtained a singular minimizer $u: \R^3 \to \R^2$ with an isolated singularity.

In the vector-valued setting, energy functionals satisfying weaker notions of convexity such as quasiconvexity and polyconvexity are also frequently considered, since they appear more naturally in many physical contexts; all these weaker notions coincide with classical convexity on scalar functions. For minimizers of strongly quasiconvex functionals, a significantly weaker bound on the singular set is currently known, see \cite{Evans86,KM07,DeFilippisKristensenMingioneStrong2026}. Meanwhile, in the recent work \cite{GT26}, Guerra \& Tione constructed a
Lipschitz singular minimizer in two dimensions for a smooth strongly
polyconvex integral.

When compared to the aforementioned vector-valued setting, the additional gradient structure of vectorial minimizers that we have herein both allows us to potentially hope for better regularity and structure than that for gradient-dependent energies on vector-valued functions, and prevents us from using the singular examples (which are not gradients) constructed in the latter setting. One can also see the difference at the level of the Euler-Lagrange equation, for which our test vector fields also have a gradient structure, therefore yielding the scalar weak fourth-order equation $\divdiv DF(D^2u)=0$, rather than the second-order system associated
with arbitrary vector-valued variations.

The construction of the energy and associated singular minimizer herein is instead motivated by partial regularity results for Hamiltonian stationary graphs. Building on earlier regularity results \cite{CW19, BW19, BCW23, Bha26, BO2025} in that setting, together with Orriols, the present authors proved a sharp
codimension-five partial regularity result for Lipschitz Hamiltonian
stationary graphs with locally $W^{1,2}$ Lagrangian angle in \cite{BOS26}. This in particular yields sharpness of the previously known partial regularity \cite{Duggan1986,Dimler2023} for Lipschitz special Lagrangian graphs. The sharpness is realized by a construction of a singular
Lipschitz graphical special Lagrangian cone in dimension five, 
using an ansatz derived from the two-homogeneous normalization of Cartan's isoparametric cubic \eqref{eq:potential}. Constructions modeled on the Cartan cubic have previously appeared in results yielding singular solutions to uniformly elliptic second order equations by Nadirashvili \& Vl\u{a}du\c{t} \cite{NV07} and Nadirashvili, Tkachev \& Vl\u{a}du\c{t} \cite{NTV12} (see also \cite{NV13}), as well as in results producing new examples of singular minimal hypercones \cite{FiresterTsiamis26} and solutions to the one-phase Bernoulli free boundary problem \cite{HW15,Wan16,FTZ26,CV2026}.

The novelty herein is to realize \eqref{eq:potential} itself as a minimizer of a suitable energy functional \eqref{eq:energy}. The freedom of tailoring the energy to this particular singular function gives us an analogue to that in \cite{NTV12}, though the variational problem considered here is different. This is somewhat in contrast to what we do in \cite{BOS26}, where instead the PDE is fixed to be the zero-phase special Lagrangian equation, which is compensated by the flexibility of seeking a solution that is a real analytic reparameterization of the level sets of the Cartan cubic foliation.

\subsection*{Overview of proof \& structure of paper}
In order to find an energy functional $F$ satisfying \eqref{eq:ellipticity} for which \eqref{eq:potential} is a minimizer, we follow the null-Lagrangian approach of \v{S}ver\'ak \& Yan \cite{SY00}. Namely, we seek a quadratic null Lagrangian $L$ (see Section \ref{positivity}) satisfying the uniform non-degeneracy condition
\begin{equation}\label{eq:secantintro}
 L(M-N)\geq\sigma|M-N|^2\qquad \text{for $M,N\in K := \{D^2u(\xi): \xi\in\Sscr_1(3)\}$}\,,
\end{equation}
for a suitable choice of $\sigma>0$, where $u$ is our candidate singular minimizer \eqref{eq:potential}. To this end, we exploit some key identities satisfied by the Cartan cubic to first find two canonical quadratic null-Lagrangians $L_0$, $L_1$ that have the same $SO(3)$ conjugation action invariance as our candidate singular minimizer. These computations are performed in Section \ref{Cartan} and the beginning of Section \ref{positivity}. Much like in \cite{MS16}, we proceed to check that the degeneracies of $L_0$ and $L_1$ compensate one another. The remainder of Section \ref{positivity} is then dedicated to identifying the degenerate directions of $L_1$ in $K$, and showing that $L_0$ is uniformly nondegenerate there, thus yielding the desired non-degeneracy \eqref{eq:secantintro} for a suitable linear combination $L = L_0 + b L_1$.

This is a sufficient condition in order to guarantee the existence of an extension to a uniformly convex integrand $F$ satisfying
\[
    F=L \quad \text{and} \quad DF = DL \qquad \text{on $K$},
\]
as already observed in \cite{SY00,MS16}. In Section~\ref{main_proofs}, we recall this extension result and verify that it applies herein and yields an extension satisfying our desired properties, therefore allowing us to conclude our main results.

The fact that $L$ is a null-Lagrangian whose gradient agrees with that of $F$ on $K$ then immediately yields the weak Euler-Lagrange equation
\[
 \divdiv DF(D^2u)=\divdiv DL(D^2u)=0
\]
for our choice of $u$ by construction.

\subsection*{Acknowledgements}
 A.B. acknowledges the support of NSF grant DMS-2350290, the Simons Foundation grant MPS-TSM-00002933, and a Bill Guthridge Fellowship from UNC-Chapel Hill. A.S.~is grateful for the generous support of Dr.~Max R\"ossler, the Walter Haefner Foundation, and the ETH Z\"urich Foundation. This research was conducted during the period A.S.~served as a Clay Research Fellow. 

\subsection*{Note on LLM Usage} The authors came up with the proof strategy and ideas underlying this paper in 2024, namely, to find an appropriate null Lagrangian that respects the same symmetries as some candidate singular minimizer, and then extend to a uniformly convex energy.
The only missing ingredient at that stage was a suitable singular candidate for a minimizer. Following the completion of the recent work \cite{BOS26}, the authors realized that the Cartan potential $P(X)/|X|$ provides precisely such a candidate, thereby completing the strategy. At this stage, AI was used as a computational aid for some of the algebraic calculations involved in identifying suitable quadratic null Lagrangians. The authors reworked these computations extensively, substantially simplifying, modifying, and entirely rewriting the arguments. The idea behind the examples in Remark \ref{r:strict-conv} was suggested by AI. All the text in this paper was written by the authors.

\section{Notation} \label{notation}
Unless otherwise specified, positive constants $c$ and $C$
depend only on the normalization of $P$ in \eqref{eq:potential}; note that there is no dependence on $n$ since we fix $n=5$. We write $\Sym(n)$ for the space of symmetric $n\times n$ matrices, while $\Sym_0(n)$ denotes the trace-free subspace of these. We let $\tr M$ denote the trace of an $n \times n$ matrix $M$. We write $\mathfrak{S}(n)$ to denote the space of self-adjoint linear operators from $\Sym_0(n)$ to $\Sym_0(n)$, which may be identified with $\Sym(\tfrac{n(n+1)}{2}-1)$. We write either $\langle X,Y\rangle=\tr(XY)$ or $X:Y$ for the Frobenius inner product on $\Sym(n)$. We let $|\cdot|$ denote the norm induced by the Frobenius inner product, and $\|\cdot \|_{\op}$ to denote the operator norm; we may equip both $\Sym_0(n)$ and $\mathfrak{S}(n)$ with both. Note that $\sqrt{\frac{n}{n-1}}\|\cdot\|_{\op}\leq |\cdot| \leq C(n) \|\cdot\|_{\op}$ on $\Sym_0(n)$. The identity operators on $\Sym_0(n)$ and $\R^n$ are denoted by $\Id_{\Sym_0(n)}$ and $\Id_n$, respectively. We let $\Sscr_1(n):= \{ X\in \Sym_0(n) : |X|=1\}$ denote the unit sphere in $\Sym_0(n)$. For $X,Y\in \Sym_0(n)$, we further use $X\otimes Y$ to denote the tensor product, which can be seen as a linear operator on $\Sym_0(n)$, given by $(X\otimes Y) Z = (X:Z) Y$.

\section{Key identities for $P$ and $u$}\label{Cartan}
We begin by recording some useful elementary identities satisfied by the Cartan cubic $P$ in \eqref{e:Cartan-cubic} and its associated two-homogeneous normalization $u$ in \eqref{eq:potential}. First of all, note that the three-homogeneity of $P$ in particular yields
\begin{equation}\label{e:P-homog}
    DP(X)\cdot X=3P(X), \qquad D^2 P(X)X =2DP(X)\,.
\end{equation}
Moreover, recall that via the map
\begin{equation*}\label{eq:map-to-Sym_0(3)}
Q(x):=
\begin{pmatrix}
2x_1&\sqrt3x_3&\sqrt3x_4\\
\sqrt3x_3&-x_1+\sqrt3x_2&\sqrt3x_5\\
\sqrt3x_4&\sqrt3x_5&-x_1-\sqrt3x_2
\end{pmatrix},
\end{equation*}
we may identify $\R^5$ with the space $\Sym_0(3)$ of symmetric trace-free matrices. Moreover, since
\begin{equation*}\label{eq:Q-invariants}
  P(x)=\frac12\det Q(x),
  \qquad
  |x|^2=\frac16\tr(Q(x)^2),
\end{equation*}
the conjugation action $\rho_R(Q):=RQR^T$, $R\in \mathrm{SO}(3)$, induces an orthogonal action on $\R^5$ under which $P$ and $u$ are invariant. In particular, we may henceforth work with isometric coordinates $X= Q(x)/\sqrt{6}$, thus identifying $P$ with a function on $\Sym_0(3)$ (not relabeled), so that $|X|^2 = |x|^2$ and $P(X) = 3\sqrt{6}\det X = \sqrt{6}\tr(X^3)$.

We begin with the following lemma, which provides us with some useful basic algebraic properties for $P$ that arise from its symmetries.
\begin{lemma}\label{lem:cubic}
For every $X,Y\in \Sym_0(3)$,
\begin{align}
 DP(X)&=3\sqrt6\left(X^2-\frac{|X|^2}{3}\Id_3\right),\label{eq:gmatrix}\\
 D^2 P(X)Y&=3\sqrt6\left(XY+YX-\frac23\tr(XY)\Id_3\right),\label{eq:Bmatrix}\\
 \Delta P(X)&=0,\qquad |DP(X)|^2=9|X|^4,\qquad |P(X)|\leq|X|^3,\label{eq:eikonal}\\
 &D^2 P(X)DP(X)=18|X|^2X,\label{eq:Bidentities}\\
 \ip{D^2 P(X)}{D^2 P(Y)}&=126\ip{X}{Y},\qquad
 |D^2 P(X)| = \sqrt{126} |X|.\label{eq:Bnorm}
\end{align}
\end{lemma}
\begin{proof}
The identities \eqref{eq:gmatrix}, \eqref{eq:Bmatrix}, and the first two identities in \eqref{eq:eikonal} all follow directly from differentiating $P(X) = 3\sqrt{6}\det X$, while the final estimate in \eqref{eq:eikonal} follows from \eqref{eq:gmatrix}, the three-homogeneity relation in \eqref{e:P-homog}, and Cauchy-Schwarz.
As for the identity \eqref{eq:Bidentities}, this is an immediate consequence of differentiating the second identity in \eqref{eq:eikonal}.

Finally, to obtain the relations in \eqref{eq:Bnorm}, we combine the elementary Hilbert space identity
\[
    \langle A, B \rangle = \frac{1}{4}\big(\|A+B\|^2 - \|A-B\|^2\big)\,,
\]
with Bochner's identity and \eqref{eq:eikonal}.

\end{proof}

Let us now use the identities and estimates in Lemma \ref{lem:cubic} in order to characterize the behavior of the Hessian of $u$. Firstly, note that by the two-homogeneity of $u$, it suffices to understand $D^2 u$ on the unit sphere in $\Sym_0(3)$ (which corresponds to $\Sbb^4$ by our isometric choice of coordinates).
Differentiating $u$, we find
\begin{align*}
 D^2 u(X) &=D^2 P(X) |X|^{-1} - DP(X)\otimes\frac{X}{|X|^3} - \frac{X}{|X|^3}\otimes DP(X)  \\
 &\qquad-P(X) |X|^{-3}\Id_{\Sym_0(3)}+3P(X) |X|^{-5} X\otimes X,\notag
\end{align*}
and thus
\begin{align}
    D^2 u(X) &=D^2 P(X) - DP(X)\otimes X - X\otimes DP(X) -P(X) \Id_{\Sym_0(3)} \label{eq:Hformula}\\
    &\qquad+3P(X) X\otimes X\notag
\end{align}
for any $X\in \Sscr_1(3):= \{ X\in \Sym_0(3) : |X|=1\}$.

Next, we record some important properties of $D^2 u$ that will be particularly important later, when extending the non-degenerate quadratic null Lagrangian that we end up with to a uniformly convex integrand.

\begin{lemma}\label{lem:Hbounds}
For $X\in\Sscr_1(3)$,
\begin{equation}\label{eq:Htrace}
 \Delta u(X)=-8P(X),\qquad |D^2 u(X)|^2=72+38P(X)^2.
\end{equation}
In particular,
\begin{equation}\label{eq:Hbounds}
 \sqrt{72}\leq|D^2 u(X)|\leq C_{Hess},
\end{equation}
where 
$C_{Hess}=\max_{\Sscr_1(3)}|D^2 u|<\infty$
only depends on $P$.
\end{lemma}
\begin{proof}
Let $X\in\Sscr_1(3)$. To see that $\Delta u =-8P$, we simply take the trace in \eqref{eq:Hformula}, and combine with the properties $\Delta P=0$ and \eqref{e:P-homog}, as well as the fact that the dimension of $\Sym_0(3)$ is 5. As for the identity for $|D^2 u(X)|^2$, we may rewrite \eqref{eq:Hformula} as $D^2 u(X)=D^2 P(X)-A_0+B_0$, with
\begin{equation}\label{e:A-and-C}
 A_0=DP(X)\otimes X+X\otimes DP(X),\qquad
 B_0=-P(X)\Id_{\Sym_0(3)}+3P(X) X\otimes X.
\end{equation}
Expanding $\langle D^2 u(X), D^2 u(X)\rangle$ via this decomposition, and invoking Lemma~\ref{lem:cubic}, \eqref{e:P-homog} and the identity $\ip{a\otimes b}{c\otimes d}=(a\cdot c)(b\cdot d)$, we obtain
\[
 |D^2 u(X)|^2=|D^2 P(X)|^2+|A_0|^2+|B_0|^2
       -2\ip{D^2 P(X)}{A_0}+2\ip{D^2 P(X)}{B_0}-2\ip{A_0}{B_0}
\]
where
\[
 \begin{gathered}
 |D^2 P(X)|^2=126,\qquad |A_0|^2=18+18P(X)^2,\qquad |B_0|^2=8P(X)^2,\\
 \ip{D^2 P(X)}{A_0}=36,\qquad \ip{D^2 P(X)}{B_0}=18 P(X)^2,\qquad \ip{A_0}{B_0}=12 P(X)^2.
 \end{gathered}
\]
This yields the second identity in \eqref{eq:Htrace} as desired. The bounds in \eqref{eq:Hbounds} follow immediately.

\end{proof}

We conclude this section with the following lemma, which guarantees smoothness and embeddedness of the set $K$ given by
\begin{equation}\label{e:K}
    K := \{D^2 u(X) : X \in \Sscr_1(3)\},
\end{equation}
cf. \eqref{eq:secantintro}. This is one of the properties that we will require to extend our final quadratic null Lagrangian to a uniformly convex integrand.

\begin{lemma}\label{lem:embedding}
The map $D^2 u :\Sscr_1(3)\to\mathfrak{S}(3)$ satisfies
\begin{equation}\label{eq:Hlipschitz}
 |D^2 u(X)-D^2 u(Y)|^2 \geq \frac{45}{2}|X-Y|^2
 \qquad X,Y\in \Sscr_1(3).
\end{equation}
In particular, $D^2 u |_{\Sscr_1(3)}$ is a smooth embedding, i.e. $K$ is a smooth embedded submanifold.
\end{lemma}
\begin{proof}
Fix $X,Y\in \Sscr_1(3)$. Let us directly compute the expansion of $|D^2 u(X)-D^2 u(Y)|^2$, analogously to the computation in the proof of Lemma \ref{lem:Hbounds}. Let $A_0$ and $B_0$ be as in \eqref{e:A-and-C}, and note that $D^2 P$ is linear. This time, it will be more convenient to keep a piece of the $D^2 P$ term with $A_0-B_0$, in order to ensure that we end up with a positive lower bound while crudely estimating the square term involving $A_0-B_0$ below by zero. More precisely, we have
\begin{align*}
 |D^2u(X)-D^2u(Y)|^2
 &= \tfrac{3}{4}|D^2P(X-Y)|^2
    \\
    &\qquad+|\tfrac{1}{2} D^2 P(X-Y)-(A_0(X)-A_0(Y))+ B_0(X)-B_0(Y)|^2\\
 &\qquad -\ip{D^2P(X-Y)}{A_0(X)-A_0(Y)} \\
 &\qquad+\ip{D^2P(X-Y)}{B_0(X)-B_0(Y)} .
\end{align*}
By Lemma \ref{lem:cubic}, the symmetry of $D^3P$ (cf. \eqref{e:3rd-der-symm}), and \eqref{e:P-homog}, we have
\begin{align*}
 |D^2P(X-Y)|^2 &=126|X-Y|^2,\\
 \ip{D^2P(X-Y)}{A_0(X)}
 &=2\ip{D^2P(X)DP(X)}{X-Y}
 =36\ip{X}{X-Y},\\
 \ip{D^2P(X-Y)}{B_0(X)}
 &=3P(X)\ip{D^2P(X)X}{X-Y}
 =6P(X)\ip{DP(X)}{X-Y},
\end{align*}
Inserting these into the above identity, and throwing away the square term involving $A_0 - B_0$, yields
\begin{equation}\label{e:Hess-diff}
 |D^2u(X)-D^2u(Y)|^2 \geq \frac{117}{2}|X-Y|^2 +6\ip{P(X)DP(X)-P(Y)DP(Y)}{X-Y}.
\end{equation}
In order to bound the final term on the right-hand side above, we proceed via the Fundamental Theorem of Calculus, using a lower bound on $D(P(X) DP(X))$. More precisely, notice that
\begin{align*}
 \ip{P(X)DP(X)-P(Y)DP(Y)}{X-Y} \\
= \int_0^1 D\big(P(Y + t(X-Y)) DP(Y+t(X-Y))\big)[X-Y,X-Y]\,dt.
\end{align*}
Moreover, note that $ D(P(Z) DP(Z)) = DP(Z) \otimes DP(Z) + P(Z) D^2 P(Z)$. Again invoking the symmetry of $D^3P$, we have
\[
 |\ip{D^2P(Z)W}{W}|
 =2|\ip{DP(W)}{Z}|
 \leq 6|Z||W|^2.
\]
Combining this with the fact that $DP(Z) \otimes DP(Z)$ is positive semidefinite and that $Z\in [X,Y] \subset \{W\in \Sym_0(3) : |W|\leq 1\}$, and exploiting the final inequality \eqref{eq:eikonal}, we therefore obtain 
\[
    D(P(Z) DP(Z)) \geq - 6 \Id_{\Sym_0(3)}.
\]
Inserting this lower bound into the integrated identity above, and combining with \eqref{e:Hess-diff}, we arrive at \eqref{eq:Hlipschitz} as desired, which in turn tells us that $D^2 u$ is an injective immersion on $\Sscr_1(3)$. The fact that it is a smooth embedding follows immediately by compactness of ${\Sscr_1(3)}$.

\end{proof}

\section{A non-degenerate quadratic null Lagrangian}\label{positivity}
We are now in a position to construct a quadratic null Lagrangian for the Hessian, i.e. a quadratic form $L$ for which the weak double-divergence form equation
\begin{equation*}\label{eq:nullidentity}
 \int DL(D^2v) : D^2 \varphi =0 \qquad \text{for every $\varphi \in C_c^\infty(B_1)$}
\end{equation*}
holds for every $v\in W^{2,2}_\loc$. See \cite{BCO81} for more details about the theory of null Lagrangians. Taking the Fourier transform (as a tempered distribution) of $\divdiv DL(D^2 (v \chi))$ for $\chi$ a smooth compactly supported function that is identically 1 on $\supp\varphi$, we immediately see that a sufficient condition for $L$ to be a null Lagrangian is thus the algebraic relation
\begin{equation}\label{eq:rankonezero}
 L(\xi \otimes \xi)=0\qquad \text{for every $\xi\in \Sym_0(3) \cong \R^5$.}
\end{equation}
Inspired by the computations in \cite[Proof of Lemma 2.1, Section 3.1]{SY00}, our construction of a suitable null Lagrangian will rely on understanding the irreducible $SO(3)$-invariant subspaces. The expression \eqref{eq:Hformula} suggests a canonical decomposition of $\mathfrak{S}(3)$ into
\begin{equation*}\label{e:subspace-decomp}
    \mathfrak{S}(3) = D^2 P(\Sym_0(3)) \oplus V_0 \oplus \langle \Id_{\Sym_0(3)} \rangle\,,
\end{equation*}
where 
\begin{align*}
    V_0 := (
D^2P(\Sym_0(3))
\oplus \langle\Id_{\Sym_0(3)}\rangle)^\perp \\
=(D^2P(\Sym_0(3)))^\perp
\cap \{M\in \mathfrak{S}(3) : \tr M = 0\},
\end{align*}
where the orthogonal complement is within $\mathfrak{S}(3)$. It is easy to verify that these $SO(3)$-invariant subspaces are irreducible. It is a classical fact (a particular case of Schur's Lemma) from invariant operator theory (see e.g. \cite{Weyl39}) that any $SO(3)$-invariant quadratic form $L$ on $\mathfrak{S}(3)$ must therefore be of the form
\begin{equation*}\label{e:L-invariant-decomp}
    L = a_1|M_1|^2 + a_2 |M_2|^2 + a_3 |M_3|^2,
\end{equation*}
for $M_1 \in D^2 P(\Sym_0(3))$, $M_2 \in V_0$, $M_3 \in \langle \Id_{\Sym_0(3)} \rangle$.
This decomposition is precisely what we will make use of in order to make appropriate choices of quadratic null Lagrangians.

Now, let us characterize $(D^2P(\Sym_0(3)))^\perp$ in a more convenient way. Let $\cT:\mathfrak{S}(3)\to \Sym_0(3)$ be the adjoint operator of $D^2 P$, namely
\begin{equation*}\label{eq:Tdefinition}
 \ip{\cT(M)}{Y}=\ip{M}{D^2 P(Y)}\qquad \text{for $M\in\mathfrak{S}(3),\ Y\in \Sym_0(3)$},
\end{equation*}
{so that 
\[
    (D^2P(\Sym_0(3)))^\perp = \{M\in \mathfrak{S}(3) : \Tcal(M) = 0\}\,.
\]
In conclusion,
\[
    V_0 = \{M\in \mathfrak{S}(3) : \tr M = 0, \ \Tcal(M)=0\}.
\]
Let us further understand how $\Tcal$ acts on $K$ in terms of $X\in\mathscr S_1(3)$; this will be crucial in order to verify the non-degeneracy conditions required from our choices of quadratic null Lagrangians.

Recalling \eqref{eq:Hformula}, this reduces to characterizing its individual action on each of the components of the Hessian therein. First of all, in light of \eqref{eq:eikonal} and \eqref{eq:Bnorm} from Lemma \ref{lem:cubic}, we have
\begin{equation}\label{eq:Tnorm}
 \cT(D^2 P(X))=126X,\qquad \|\cT\|_{\op}=\sqrt{126}, \qquad \cT(\Id_{\Sym_0(3)})=0.
\end{equation}
As for the action of $\Tcal$ on rank-one tensors, the symmetry of $D^3P$ as a constant 3-tensor yields
\begin{equation}\label{e:3rd-der-symm}
 \langle \cT(X\otimes X), Y \rangle = \langle D^2 P(Y) X, X \rangle = D^3 P(0)[Y,X,X] = \langle D^2 P(X) X, Y \rangle,
\end{equation}
and thus, combining with the homogeneity \eqref{e:P-homog} of $D^2 P$ we have
\begin{equation}\label{eq:Trankone}
    \cT(X\otimes X) = 2DP(X).
\end{equation}}
Finally, for $X \in \Sscr_1(3)$, \eqref{eq:Bidentities} additionally yields
\begin{equation}\label{e:Tsymm-piece}
 \cT\bigl(DP(X)\otimes X+X\otimes DP(X)\bigr)=36X.
\end{equation}
Combining \eqref{eq:Tnorm}, \eqref{eq:Trankone} and \eqref{e:Tsymm-piece}, we therefore obtain
\begin{equation}\label{eq:TH}
 \cT(D^2 u(X))=90X+6P(X)DP(X).
\end{equation}

Let us now use the above decomposition into $SO(3)$-irreducible invariant subspaces in order to find a pair of null Lagrangians whose degeneracies compensate one another. 
We begin by identifying two quadratic null Lagrangians arising naturally from the computations in the preceding section. Consider two quadratic forms, each of which depends only on two components $M_i$ in a decomposition $(M_1,M_2,M_3)\in  D^2 P(\Sym_0(3)) \oplus V_0 \oplus \langle \Id_{\Sym_0(3)} \rangle$, for instance
\begin{align*}
 L_1(M)&=a_1|M_1|^2 + a_3 |M_3|^2,\\
 L_0(M)&= b_2|M_2|^2+ b_3|M_3|^2.
\end{align*}
These are convenient ans\"{a}tze since we avoid having to work with linear combinations of all 3 of the $M_i$, while retaining enough flexibility to solve for the coefficients in order to satisfy the null Lagrangian condition \eqref{eq:rankonezero}. However, many other choices would also work, subject to satisfying the latter constraint. Moreover, note that we only expect to find two linearly independent quadratic null Lagrangians with $SO(3)$-invariance, since there are three mutually orthogonal subspaces $D^2 P(\Sym_0(3))$, $V_0$, $\langle \Id_{\Sym_0(3)} \rangle$ from which to take linear combinations, but the null Lagrangian constraint removes one degree of freedom.

First of all, letting $X\in \Sym_0(3)$ be such that $M_1 = D^2 P(X)$ and noting that $\Tcal \equiv 0$ on $V_0 \oplus \langle \Id_{\Sym_0(3)} \rangle$, \eqref{eq:Tnorm} yields $126 X = \Tcal(M)$, and thus thanks to \eqref{eq:Bnorm} we obtain 
\[
    |M_1|^2 = |D^2 P(X)|^2 = 126|X|^2 = \frac{1}{126}|\Tcal(M)|^2.
\]
Secondly, note that \eqref{eq:eikonal} guarantees $\tr M_1 = 0$, and we have $\tr M_2 =0$ by construction. Thus, we must have $\tr M_3 = \tr M$, and since $\dim \Sym_0(3) = 5$, we obtain 
\[
    |M_3|^2 = \frac{1}{5} (\tr M)^2.
\]
Orthogonality then allows us to recover $|M_2|^2$ via
\[
    |M_2|^2 = |M|^2 - |M_1|^2 - |M_3|^2.
\]
Taking $M=\xi\otimes \xi$ and combining the above properties with \eqref{eq:Trankone}, \eqref{eq:eikonal} and the null Lagrangian constraint \eqref{eq:rankonezero} which requires $L_0$ and $L_1$ to vanish on rank-one tensors, we obtain the relations
\[
    \frac{2}{7} a_1 + \frac{1}{5} a_3 = 0, \qquad \left(1-\frac{2}{7}-\frac{1}{5}\right) b_2 + \frac{1}{5} b_3 = 0\,.
\]
Taking $a_1=1$ and $b_2=-1$ yields
\begin{align*}
 L_1(M)&=\frac1{126}|\cT(M)|^2-\frac27(\tr M)^2,\\
 L_0(M)&=\frac1{126}|\cT(M)|^2+\frac57(\tr M)^2-|M|^2.
\end{align*}
Again, this is simply a choice that we fix, for concreteness; any other one would work equally well.

By construction, we have therefore obtained the following.
\begin{proposition}\label{prop:null}
Both $L_0$ and $L_1$ are
quadratic null Lagrangians for the Hessian. They are both invariant under the action
$M\mapsto\rho_R M\rho_R^*$ induced by $R\in SO(3)$.
\end{proposition}

We will see that although $L_1$ is strictly
positive on every direction $M-N$ between distinct $M,N$ in the image $K$ of $\Sscr_1(3)$ under $D^2 u$, this positivity is not uniform and indeed
degenerates along certain infinitesimal directions as $M$ and $N$ approach one another. However, we will verify that
$L_0$ remains uniformly non-degenerate on these directions. This complementary behavior is what allows us to find a linear combination of $L_0$ and $L_1$ that is uniformly non-degenerate on $K$, as stated in the following key proposition; this is what we will need in order to extend this linearization to our nonlinear uniformly convex energy.

\begin{proposition}\label{prop:secant}
There are constants $b\geq0$ and $\sigma>0$ 
such that the quadratic null Lagrangian $L=L_0+bL_1$ satisfies
\begin{equation}\label{eq:secant}
 L(M-N)\geq\sigma|M-N|^2\qquad \text{for $M,N\in K$,}
\end{equation}
where $K$ is as in \eqref{e:K}.
\end{proposition}

The remainder of this section will be dedicated to the proof of Proposition \ref{prop:secant}, which we subdivide into a number of intermediate lemmas. In the first of these, we show strict positivity of $L_1$ on $K$.

\begin{lemma}\label{lem:strictsecants}
For distinct $X,Y\in\Sscr_1(3)$,
\begin{equation*}\label{eq:Jstrict}
 L_1(D^2 u(X)-D^2 u(Y))=\frac{1}{126}|\cT(D^2 u(X))-\cT(D^2 u(Y))|^2-8^2\cdot\frac{2}{7}|P(X)-P(Y)|^2 > 0.
\end{equation*}
\end{lemma}
\begin{proof}
Recalling \eqref{eq:TH}
and \eqref{eq:gmatrix}, for $X\in \Sscr_1(3)$ we have
\begin{equation}\label{e:TH-rewritten}
 \Tcal(D^2 u(X))=90X+18\sqrt6 P(X)\left(X^2-\frac{1}{3}\Id_3\right).
\end{equation}
Thus, if we choose an orthonormal basis $\{e_i\}_{i=1}^3$ of eigenvectors for $X$, i.e. $Xe_i = \lambda_{X,i} e_i$, we in turn have $X^2 e_i = \lambda_{X,i}^2 e_i$ and so $\{e_i\}$ is also an eigenbasis for $\Tcal(D^2 u(X))$. Letting $\{\mu_{X,i}\}$ denote the associated eigenvalues for $\Tcal(D^2 u(X))$ observe that the Rayleigh quotient characterization of the eigenvalues yields
\begin{equation}\label{e:Rayleigh}
    |\mu_{X,i} - \mu_{Y,i}| \leq \|\Tcal(D^2 u(X)) - \Tcal(D^2 u(Y))\|_{\op}.
\end{equation}
We would like to find a lower bound on $|\mu_{X,i} - \mu_{Y,i}|$ in terms of $X$ directly, rather than via the adjoint $\Tcal$. List the eigenvalues $\{\lambda_{X,i}\}_i$ monotonically increasing in size, i.e. $\lambda_{X,1} \leq \lambda_{X,2} \leq \lambda_{X,3}$. We will proceed to express everything in terms of the \emph{middle} one $\lambda_{X,2}$, though it could be done with any fixed one of them; the reason for this choice will become apparent later. Arrange $X,Y\in \Sscr_1(3)$ so that $\lambda_{X,2} \geq \lambda_{Y,2}$ without loss of generality. Now, recall that 
\[
    P(X) = 3\sqrt{6} \det X = 3\sqrt{6} \prod_{i=1}^3\lambda_{X,i}.
\]
Since $\tr X = \sum_{i=1}^3 \lambda_{X,i} =0$ and $|X|^2 = \sum_{i=1}^3 \lambda_{X,i}^2 =1$, we obtain $\lambda_{X,1}\lambda_{X,3} = \lambda_{X,2}^2 -\frac{1}{2}$, and therefore
\[
    P(X) = 3\sqrt{6} \lambda_{X,2}\left(\lambda_{X,2}^2 -\frac{1}{2}\right) =: p(\lambda_{X,2}).
\]
Inserting this into the identity \eqref{e:TH-rewritten} and taking the $e_2$ component of the latter yields
\begin{align*}
    |\mu_{X,2} - \mu_{Y,2}| &= \left|90(\lambda_{X,2}- \lambda_{Y,2}) + 18\sqrt{6}\left[p(\lambda_{X,2})\left(\lambda_{X,2}^2 -\frac{1}{3}\right) - p(\lambda_{Y,2})\left(\lambda_{Y,2}^2 -\frac{1}{3}\right)\right]\right|\\
    &=\left|90(\lambda_{X,2}- \lambda_{Y,2}) + 18\sqrt{6} \int_{\lambda_{Y,2}}^{\lambda_{X,2}} \frac{d}{ds}\left[p(s) \left(s^2 -\frac{1}{3}\right)\right]ds\right|.
\end{align*}
Let us now compute the derivative inside the integral. Since $p(s) = 3\sqrt{6} s\left(s^2 -\frac{1}{2}\right)$, we have
\begin{align*}
    18\sqrt{6}\frac{d}{ds}\left[p(s)\left(s^2 - \frac{1}{3}\right)\right] = 18^2\left(5s^4-\frac{5}{2}s^2+\frac{1}{6}\right) \\
    = \frac{d}{ds}\left[-18\sqrt{6}\cdot\frac{8}{9}p(s)-90 s\right]
  +54 s^2(1+30s^2),
\end{align*}
and thus
\begin{equation}\label{e:mu-eval-lb}
     |\mu_{X,2} - \mu_{Y,2}| = \left|-16\sqrt{6}[p(\lambda_{X,2}) - p(\lambda_{Y,2})] + 54 \int_{\lambda_{Y,2}}^{\lambda_{X,2}} s^2(1+30s^2)\, ds\right|.
\end{equation}
Now our choice of the middle eigenvalue $\lambda_{X,2}$ comes into play. We claim that for this one, $p'(s)\leq 0$ on the segment $[\lambda_{Y,2},\lambda_{X,2}]$ and so in particular $-[p(\lambda_{X,2}) - p(\lambda_{Y,2})]\geq 0$.
Indeed, notice that $p'(s)=3\sqrt6\left(3s^2-\frac{1}{2}\right)$, so
\[
    p'(\lambda_{X,2})=3\sqrt6\left(3\lambda_{X,2}^2-\frac{1}{2}\right) = 3\sqrt{6}(\lambda_{X,2} -\lambda_{X,1})(\lambda_{X,2}-\lambda_{X,3}) \leq 0,
\]
where the second identity is again exploiting the relations coming from $\tr X =0$ and $|X|^2 =1$. This in turn implies that $\lambda_{X,2} \in [-\frac{1}{\sqrt{6}},\frac{1}{\sqrt{6}}]$, and so the claim follows.

Inserting this into \eqref{e:mu-eval-lb} and combining with \eqref{e:Rayleigh}, we deduce that
\begin{align*}
      |\Tcal(D^2 u(X)) - \Tcal(D^2 u(Y))|^2 \geq \frac{3}{2}\|\Tcal(D^2 u(X)) - \Tcal(D^2 u(Y))\|^2_{\op} \\
      \geq 9\cdot 16^2 (P(X)-P(Y))^2\,.
\end{align*}
Combining this with the Laplacian identity in Lemma \ref{lem:Hbounds}, the conclusion follows.

\end{proof}

Note that the positivity in Lemma \ref{lem:strictsecants} is not uniform, and might degenerate as $X,Y$ approach one another. Let us therefore identify these infinitesimal directions on which $L_1$ degenerates and show that on these, $L_0$ is bounded uniformly below.

\begin{lemma}\label{lem:tangent}
For $X\in\Sscr_1(3)$ and $v\in T_X\Sscr_1(3)$,
\[
 L_1(D(D^2 u)(X)[v]) \equiv L_1((D^3 u)(X)v)\geq0,
\]
If $v\neq0$ and equality holds, then
\begin{equation*}\label{eq:tangentratio}
 \frac{L_0(D^3 u(X)v)}{|D^3 u(X)v|^2}=c_0>0,
\end{equation*}
where $c_0$ is a universal constant depending on $P$. 
Here, $(D^3 u)(X)v = (D^3 u)(X)[v,\cdot,\cdot]$ is an element of $\mathfrak{S}(3)$.
\end{lemma}
\begin{proof}

    Fix $X\in\Sscr_1(3)$ and $v\in T_X\Sscr_1(3)$. We will proceed to directly compute the behavior of $L_1$ on $D^3 u(X)v$. First of all, recall from definition that we have $L_1(D^3 u(X) v) = \frac{1}{126} |\Tcal(D^3 u(X) v)|^2 - \frac{2}{7}\big(\tr[D^3 u(X) v]\big)^2$. Using linearity of the trace operator and differentiating the Laplacian identity in Lemma \ref{lem:Hbounds} yields
    \begin{equation}\label{e:tr-D^3u}
        \tr [D^3 u(X) v] = - 8 DP(X) v,
    \end{equation}
    and thus
    \begin{equation}\label{e:L1-D^3u}
        L_1(D^3 u(X) v) = \frac{1}{126} |\Tcal(D^3 u(X) v)|^2 - \frac{128}{7}\big(DP(X) v\big)^2
    \end{equation}
    
    Noting that $\Tcal$ is a linear operator, we may simply differentiate the identity \eqref{eq:TH} in order to analyze the behavior of $\Tcal(D^3 u(X)v)$. This yields
    \begin{equation}\label{e:T-D^3u}
        \Tcal(D^3 u(X)v) = 90 v + 6 (DP(X) v) DP(X) + 6 P(X) D^2 P(X) v.
    \end{equation}
    Keeping in mind that we wish to realize \eqref{e:L1-D^3u} as a non-negative expression, let us try to complete the square. Combining \eqref{e:T-D^3u} with Lemma \ref{lem:cubic} and the observation that $\langle X, v \rangle=0$, we obtain
    \begin{align*}
        \langle \Tcal(D^3 u(X) v), DP(X) \rangle &= 90 DP(X) v + 6 |DP(X)|^2 DP(X) v + 6 P(X) D^2 P(X) DP(X) v \\
        &= 90 DP(X) v + 54 DP(X) v + 108 P(X) \langle X, v \rangle \\
        &=144 DP(X) v.
    \end{align*}
    Thus,
    \begin{align*}
        L_1((D^3 u)(X)v) &= \frac{1}{126}\left|\Tcal(D^3 u(X)v) - 16 (DP(X)v) DP(X) \right|^2 \\
        &= \frac{1}{126}\left|90v + 6 P(X) D^2 P(X) v - 10 (DP(X) v) DP(X)\right|^2 \\
        & \geq 0
    \end{align*}
    as desired.

    Let us now investigate when equality holds. Taking an inner product of the expression inside the square with $X$ and again invoking Lemma \ref{lem:cubic} and $\langle X, v \rangle=0$, a necessary condition for equality to hold is $P(X) D P(X) v = 0$. Therefore, either $P(X)=0$ or $D P(X) v = 0$. In the latter case, inserting this condition back into the squared expression above yields $|90v + 6 P(X) D^2 P(X) v|^2=0$, which in turn forces $v=0$ thanks to the bounds on $D^2 P(X)$ and $P(X)$ from Lemma \ref{lem:cubic}.

    It therefore remains to investigate the behavior of $\frac{L_0(D^3 u(X)v)}{|D^3 u(X)v|^2}$ when $P(X)=0$ and $v\neq 0$. In this case, the expression in the square tells us that $v = a DP(X)$ for some $a\neq0$. Without loss of generality, we may take $a=1/3$ so that $|v|=1$. Once again using Lemma \ref{lem:cubic}, this choice further yields $DP(X)v=3$ and
	$D^2P(X)v=6X$. Now, since 
    \begin{equation}\label{e:L1-L0-rel}
        L_0(D^3 u(X)v) = L_1(D^3 u(X)v) + (\tr[D^3 u(X)v])^2  - |D^3 u(X)v|^2
    \end{equation}
    and we know $L_1(D^3 u(X)v)=0$ and have already computed $(\tr[D^3 u(X)v])^2$ in \eqref{e:tr-D^3u}, the only thing left is to compute $|D^3 u(X)v|^2$.
    
    Differentiating \eqref{eq:Hformula} in the direction $v$ and substituting in these identities together with $P(X)=0$ yields
	\begin{align*}
	D^3u(X)v &= D^3P(X)v - 3 X\otimes X-\frac{2}{3} DP(X)\otimes DP(X) -3\Id_{\Sym_0(3)},
	\end{align*}
    and thus, combining with the fact that $D^3P(X)v = D^2P(v)$ since $P$ is a cubic polynomial, we obtain
    \begin{align*}
	   |D^3u(X)v|^2 &= |D^2P(v)|^2 +9|X|^4 +\frac{4}{9}|DP(X)|^4 + 9|\Id_{\Sym_0(3)}|^2 \\
       &-6\langle D^3P(X)v,X\otimes X\rangle
        -\frac43\langle D^3 P(X)v,DP(X)\otimes DP(X)\rangle -6\Delta P(v) \\
    &+4\langle X\otimes X,DP(X)\otimes DP(X)\rangle +18|X|^2 +4 |DP(X)|^2.
	\end{align*}
    Now, recalling the tensor identity $\ip{a\otimes b}{c\otimes d}=(a\cdot c)(b\cdot d)$ and differentiating \eqref{eq:Bidentities} in the direction $v$, we obtain
    \begin{align*}
        \langle D^3 P(X)v,DP(X)\otimes DP(X)\rangle + \langle D^2 P(X) (D^2P(X)v), DP(X) \rangle\\
        = 36 (X \cdot v)(DP(X)\cdot X) + 18 |X|^2 DP(X) v.
    \end{align*}
    Combining this with the symmetry of $D^3 P(X)$, the properties in Lemma \ref{lem:cubic}, the homogeneity relations \eqref{e:P-homog}, the fact that $\dim \Sym_0(3) =5$, and $|X|=|v|=1$, we arrive at
    \begin{align*}
	   |D^3u(X)v|^2 = 306.
	\end{align*}
Combining this with \eqref{e:L1-L0-rel} and \eqref{e:tr-D^3u}, we conclude that
\[
 \frac{L_0(D^3u(X)v)}{|D^3u(X)v|^2}
 =\frac{270}{306}=\frac{15}{17}=:c_0>0.
\]

\end{proof}

With the competing non-degeneracies at hand, we are now in a position to conclude the proof of Proposition \ref{prop:secant}.
\begin{proof}[Proof of Proposition \ref{prop:secant}]

Our goal is to find a suitable linear combination $L=L_0 + bL_1$ satisfying \eqref{eq:secant}. Recall from Lemma \ref{lem:strictsecants} that $L_1(M-N) \geq 0$ for every $M,N\in K$; it remains to compensate this with a piece of $L_0$ to obtain uniform nondegeneracy. Fix $\sigma \in (0,c_0)$, where $c_0$ is as in Lemma \ref{lem:tangent}, and set
\[
 U=\{(M,N)\in K\times K:L_0(M-N)\geq \sigma|M-N|^2\}.
\]
By definition, on $U$ the desired lower bound holds for any choice of $b\geq 0$, and so it suffices to restrict our attention to $(K\times K)\setminus U$. We claim that
\begin{equation}\label{e:L1-nondeg-final}
 \eta:=\inf_{\substack{(M,N)\in(K\times K)\setminus U\\M\neq N}}
 \frac{L_1(M-N)}{|M-N|^2}>0.
\end{equation}
Suppose for a contradiction that $\eta=0$. Then there exists a sequence $(M_j,N_j) \in K\times K$ for which
the ratios in \eqref{e:L1-nondeg-final} converge to zero. Up to extracting a subsequence, we must therefore have $M_j, N_j \to M \in K$. By Lemma \ref{lem:embedding}, the corresponding points $X_j,Y_j \in \Sscr_1(3)$ that respectively satisfy $D^2 u(X_j) = M_j$, $D^2 u(Y_j) = N_j$ also converge (subsequentially) to the same element $X\in\Sscr_1(3)$. This induces a tangent direction $v \in T_X\Sscr_1(3)$ as follows:
\[
 \frac{X_j-Y_j}{|X_j-Y_j|}\longrightarrow v,
 \qquad |v|=1,\qquad v\in T_X\Sscr_1(3),
\]
up to yet another subsequence. Taylor expanding between $Y_j$ and $X_j$ yields
\[
 \frac{M_j-N_j}{|X_j-Y_j|} = \frac{D^2 u(X_j)- D^2 u(Y_j)}{|X_j-Y_j|}\longrightarrow D^3u(X)v.
\]
Again invoking Lemma \ref{lem:embedding}, we observe that the right-hand side is non-zero. However, by construction of our contradiction sequences, we have 
\[
    L_1(D^3u(X)v)=0, \qquad L_0(D^3u(X)v)\leq \sigma|D^3u(X)v|^2,
\]
which contradicts Lemma~\ref{lem:tangent}, since $\sigma<c_0$.
Thus we indeed must have $\eta>0$, as claimed.

Next, noticing that $L_0(M-N)\geq-|M-N|^2$ by definition, we obtain
\[
 (L_0+bL_1)(M-N)\geq(-1+b\eta)|M-N|^2
\]
for $(M,N)\in(K\times K)\setminus U$ with $M\neq N$.  
Thus, letting $b\geq(1+\sigma)/\eta$, the desired nondegeneracy \eqref{eq:secant} follows for this choice of $b$.

\end{proof}

\section{Proof of main results} \label{main_proofs}
{The final step required to conclude our main results is an extension of $L$ from $K$ to a uniformly
convex integrand $F$ whose linearization is $L$ on $K$. This will be a consequence of the following lemma, which is a special case of \cite[Lemma 2.1]{MS16} where the function $G$ therein is taken to be quadratic. We refer the reader to the proof therein, and merely clarify why the extension constructed there also yields a uniformly bounded Hessian for $F$, which is what we are assuming herein.}

\begin{lemma}\label{lem:extension}
Let $K$ be a compact smooth embedded submanifold of $\R^q$. Suppose that $L:\R^q\to\R$ is a quadratic form satisfying
\[
 L(M-N)\geq\sigma|M-N|^2\qquad M,N \in K,
\]
for some $\sigma>0$. Then there exists $F\in C^\infty(\R^q)$ satisfying the bounds \eqref{eq:ellipticity} with $\lambda = \sigma$ and $\Lambda=\Lambda(q,K,L,\sigma)<\infty$, with the property that on $K$ both $F=L$ and $DF=DL$ hold.
\end{lemma}
\begin{proof}
If $K$ is nonempty, we may 
apply \cite[Lemma 2.1]{MS16} with $G= L$, $v= DL$, and $\gamma=\sigma$ in the notation therein. This in addition yields the Hessian lower bound. Meanwhile, to obtain the uniform upper bound on the Hessian, we consider the function $H_0$ from Step 2 of \cite[Section 4.1]{MS16}, given by 
\[
    H_{0}(M)=\sup_{N \in K_\delta} \{F(N) + \langle DF(N), M-N \rangle +\frac{3}{4}\sigma |M-N|^2 \},
\]
where $K_\delta$ is a $\delta$-neighborhood of $K$, for $\delta>0$ sufficiently small, and $F$ is the local extension obtained in \cite[Proof of Lemma 2.1, Step 1]{MS16}. Expanding the square in the right-hand side above, this may be rewritten as
\[
    H_{0}(M) - \frac{3}{4} \sigma |M|^2 =\sup_{N \in K_\delta} \left\{F(N) - \langle DF(N), N \rangle + \left\langle DF(N) - \frac{3}{2} \sigma N , M \right\rangle +\frac{3}{4}\sigma |N|^2 \right\}.
\]
Noting that $DF(N)-\frac{3}{2}\sigma N$ is uniformly bounded over $N\in K_\delta$, we deduce that $H_0-\frac{3}{4}\sigma|\cdot|^2$ is globally
Lipschitz, with Lipschitz constant at most
$\sup_{N\in K_\delta}|DF(N)-\frac32\sigma N|$. 

Now, for fixed $\eps>0$ suitably small, let $H_\eps := H_0 * \rho_\eps$ for a standard mollifier $\rho_\eps$ be as in \cite[Proof of Lemma 2.1, Step 2]{MS16}. Note that we may ensure $\|D\rho_\varepsilon\|_{L^1} \leq C \eps^{-1}$. Taking the Hessian of both sides above, observe that this mollification satisfies the uniform Hessian bounds
\[
 \sup_ {M\in  \mathbb R^{q}}\|D^2 H_\eps(M)\|_{\op}
 \leq \frac32\sigma
 +C\eps^{-1}\mathrm{Lip}\Big(H_0 - \frac{3}{4} \sigma |M|^2\Big) <\infty.
\]
 The final extension is then an interpolation between $H_0$ on $K_\delta$ and $H_\eps$ on $\R^q \setminus K_{2\delta}$, so the above estimates indeed guarantee a uniform upper Hessian bound.
\end{proof}

Let us now collect everything together to conclude our main results.

\begin{proof}[Proof of Theorem \ref{thm:main}]
Let $L$ be the null Lagrangian from Proposition \ref{prop:secant}. Since Lemma \ref{lem:embedding} in addition tells us that the set $K$ is a smooth, compact embedded submanifold of $\mathfrak{S}(3)\cong \Sym(5)$, we may apply Lemma \ref{lem:extension} to obtain a smooth extension $F$ that satisfies \eqref{eq:ellipticity} and $F=L$, $DF=DL$ on $K$. In particular, by construction we have
\begin{equation}\label{e:null-Lag-L-F}
 \int_{B_1}DF(D^2u):D^2\varphi\dd x
 =\int_{B_1}DL(D^2u):D^2\varphi\dd x=0 
 \qquad \text{for every $\varphi\in C_c^\infty(B_1)$.}
\end{equation}

Since $DF(D^2u)$ is uniformly bounded, \eqref{e:null-Lag-L-F} remains valid for $\varphi\in W^{2,2}_0(B_1)$ by density. Using Taylor's Theorem with integral remainder together with the uniform convexity of $F$ (i.e. the lower bound in \eqref{eq:ellipticity}), we obtain
\begin{align*}
   & F(D^2u+D^2\varphi)-F(D^2u) \\
    &= DF(D^2u):D^2\varphi + \int_0^1 (1-t) D^2 F(D^2 u+t D^2 \varphi)[D^2 \varphi,D^2 \varphi]\dd t \\
    &\geq DF(D^2u):D^2\varphi+\frac\lambda2|D^2\varphi|^2,
\end{align*}
from which \eqref{eq:strictminimum} follows immediately after integration over $B_1$.
The uniqueness of $u$ as a minimizer follows immediately from the uniform convexity \eqref{eq:strictminimum}. Note that the stated boundary data for $u$ is due to its two-homogeneity.
\end{proof}

\begin{proof}[Proof of Corollary \ref{cor:dimensions}]
This follows immediately by splitting the variables in $\R^n$ as $x=(x',y)\in\R^5\times\R^{n-5}$ and considering the candidate minimizer $x \mapsto u(x')$ as a function on $\R^n$. This function, which we still call $u$ at the risk of abusing notation, clearly satisfies \eqref{e:higher-dim-spine} for its singular set, by construction. The associated uniformly convex integrand can be taken to be 
\[
    \Sym(n) \ni M \mapsto F(A)+\frac{1}{2}(|M|^2-|A|^2),
\]
where $F$ is the integrand on $\Sym(5)$ from Theorem \ref{thm:main} and $A$ is the upper left $5\times5$ block of $M$. It is immediate to see that this satisfies \eqref{eq:ellipticity} with constants $\min\{ \lambda, 1\}$ and $\max\{ \Lambda, 1\}$, where $\lambda, \Lambda$ are the constants associated to $F$. It is easy to verify that this integrand satisfies the desired properties; we leave the details to the reader.
\end{proof}

\bibliographystyle{amsalpha}
\bibliography{library}

@misc{BO2025,
  author  = {Bhattacharya, Arunima and Ogden, W. Jacob},
  title   = {Optimal regularity for {H}\"older continuous {H}amiltonian stationary {L}agrangian graphs},
  year = {2025},
  note = {Preprint, arXiv:2510.22756},
  eprint = {2510.22756}, archivePrefix = {arXiv}, primaryClass = {math.AP},
  url = {https://arxiv.org/abs/2510.22756}
}

@article{BCO81,
  author = {Ball, John M. and Currie, J. C. and Olver, Peter J.},
  title = {Null {Lagrangians}, weak continuity, and variational problems of arbitrary order},
  journal = {J. Funct. Anal.},
  volume = {41}, number = {2}, pages = {135--174}, year = {1981},
  doi = {10.1016/0022-1236(81)90085-9}
}

@article{Bha22,
  author = {Bhattacharya, Arunima},
  title = {Regularity for critical points of convex functionals on {Hessian} spaces},
  journal = {Proc. Amer. Math. Soc.},
  volume = {150}, number = {12}, pages = {5217--5230}, year = {2022},
  eprint = {2108.01319}, archivePrefix = {arXiv}, primaryClass = {math.AP}
}

@article{Bha26,
  author = {Bhattacharya, Arunima},
  title = {Regularity for {Hamiltonian} stationary equations in {$\mathbb{R}_{n\leq4}^{n}$}},
  journal = {Math. Z.},
  volume = {312}, pages = {Paper No. 11}, year = {2026},
  doi = {10.1007/s00209-025-03907-3},
  eprint = {2504.18013}, archivePrefix = {arXiv}, primaryClass = {math.AP}
}

@article{BCW23,
  author = {Bhattacharya, Arunima and Chen, Jingyi and Warren, Micah},
  title = {Regularity of {Hamiltonian} stationary equations in symplectic manifolds},
  journal = {Adv. Math.},
  volume = {424}, pages = {Paper No. 109059}, year = {2023},
  doi = {10.1016/j.aim.2023.109059},
  eprint = {2108.00325}, archivePrefix = {arXiv}, primaryClass = {math.AP}
}

@misc{BOS26,
  author = {Bhattacharya, Arunima and Orriols, Gerard and Skorobogatova, Anna},
  title = {Sharp partial regularity of {Hamiltonian} stationary and special {Lagrangian} graphs},
  year = {2026},
  note = {Preprint, arXiv:2608.27324},
  eprint = {2608.27324}, archivePrefix = {arXiv}, primaryClass = {math.DG}
}

@article{BS25,
  author = {Bhattacharya, Arunima and Skorobogatova, Anna},
  title = {Variational integrals on {Hessian} spaces: partial regularity for critical points},
  journal = {Nonlinear Anal.},
  volume = {255}, pages = {Paper No. 113760}, year = {2025},
  doi = {10.1016/j.na.2025.113760},
  eprint = {2307.01191}, archivePrefix = {arXiv}, primaryClass = {math.AP}
}

@article{BW21,
  author = {Bhattacharya, Arunima and Warren, Micah},
  title = {Regularity bootstrapping for 4th-order nonlinear elliptic equations},
  journal = {Int. Math. Res. Not. IMRN},
  volume = {2021}, number = {6}, pages = {4324--4348}, year = {2021},
  doi = {10.1093/imrn/rny301},
  eprint = {1706.05501}, archivePrefix = {arXiv}, primaryClass = {math.AP}
}

@article{Car39,
  author = {Cartan, {\'E}lie},
  title = {Sur des familles remarquables d'hypersurfaces isoparam{\'e}triques dans les espaces sph{\'e}riques},
  journal = {Math. Z.},
  volume = {45}, pages = {335--367}, year = {1939},
  doi = {10.1007/BF01580289}
}

@article{CW19,
  author = {Chen, Jingyi and Warren, Micah},
  title = {On the regularity of {Hamiltonian} stationary {Lagrangian} submanifolds},
  journal = {Adv. Math.},
  volume = {343}, pages = {316--352}, year = {2019},
  eprint = {1611.02641}, archivePrefix = {arXiv}, primaryClass = {math.AP}
}

@article{DG57,
  author = {De Giorgi, Ennio},
  title = {Sulla differenziabilit{\`a} e l'analiticit{\`a} delle estremali degli integrali multipli regolari},
  journal = {Mem. Accad. Sci. Torino Cl. Sci. Fis. Mat. Nat. (3)},
  volume = {3}, pages = {25--43}, year = {1957}
}

@article{Evans86,
  author = {Evans, Lawrence C.},
  title = {Quasiconvexity and partial regularity in the calculus of variations},
  journal = {Arch. Rational Mech. Anal.},
  volume = {95}, number = {3}, pages = {227--252}, year = {1986},
  doi = {10.1007/BF00251360}
}

@misc{GT26,
  author = {Guerra, Andr{\'e} and Tione, Riccardo},
  title = {A singular minimizer of a smooth strongly
           polyconvex functional in two dimensions},
  year = {2026},
  note = {Preprint, arXiv:2609.32404},
  eprint = {2609.32404},
  archivePrefix = {arXiv},
  primaryClass = {math.AP},
  url = {https://arxiv.org/abs/2609.32404}
}

@article{HLN96,
  author = {Hao, Wenge and Leonardi, Salvatore and Ne{\v{c}}as, Jind{\v{r}}ich},
  title = {An example of irregular solution to a nonlinear {Euler--Lagrange} elliptic system with real analytic coefficients},
  journal = {Ann. Scuola Norm. Sup. Pisa Cl. Sci. (4)},
  volume = {23}, number = {1}, pages = {57--67}, year = {1996},
  url = {https://www.numdam.org/item/ASNSP_1996_4_23_1_57_0/}
}

@article{KM07,
  author = {Kristensen, Jan and Mingione, Giuseppe},
  title = {The singular set of {Lipschitzian} minima of multiple integrals},
  journal = {Arch. Rational Mech. Anal.},
  volume = {184}, number = {2}, pages = {341--369}, year = {2007},
  doi = {10.1007/s00205-006-0036-2}
}

@article{MS16,
  author = {Mooney, Connor and Savin, Ovidiu},
  title = {Some singular minimizers in low dimensions in the calculus of variations},
  journal = {Arch. Rational Mech. Anal.},
  volume = {221}, number = {1}, pages = {1--22}, year = {2016},
  doi = {10.1007/s00205-015-0955-x},
  eprint = {1503.00671}, archivePrefix = {arXiv}, primaryClass = {math.AP}
}

@article{NTV12,
  author = {Nadirashvili, Nikolai and Tkachev, Vladimir and Vl{\u{a}}du{\c{t}}, Serge},
  title = {A non-classical solution to a {Hessian} equation from {Cartan} isoparametric cubic},
  journal = {Adv. Math.},
  volume = {231}, number = {3--4}, pages = {1589--1597}, year = {2012},
  doi = {10.1016/j.aim.2012.07.005},
  eprint = {1111.0329}, archivePrefix = {arXiv}, primaryClass = {math.AP}
}

@article{NV07,
  author = {Nadirashvili, Nikolai and Vl{\u{a}}du{\c{t}}, Serge},
  title = {Nonclassical solutions of fully nonlinear elliptic equations},
  journal = {Geom. Funct. Anal.},
  volume = {17}, number = {4}, pages = {1283--1296}, year = {2007},
  doi = {10.1007/s00039-007-0626-7}
}

@article{NV13,
  author = {Nadirashvili, Nikolai and Vl{\u{a}}du{\c{t}}, Serge},
  title = {Singular solutions of {Hessian} elliptic equations in five dimensions},
  journal = {J. Math. Pures Appl. (9)},
  volume = {100}, number = {6}, pages = {769--784}, year = {2013},
  doi = {10.1016/j.matpur.2013.03.001},
  eprint = {1211.5812}, archivePrefix = {arXiv}, primaryClass = {math.AP}
}

@article{Nash58,
  author = {Nash, John},
  title = {Continuity of solutions of parabolic and elliptic equations},
  journal = {Amer. J. Math.},
  volume = {80}, number = {4}, pages = {931--954}, year = {1958},
  doi = {10.2307/2372849}
}

@incollection{Nec77,
  author = {Ne{\v{c}}as, Jind{\v{r}}ich},
  title = {Example of an irregular solution to a nonlinear elliptic system with analytic coefficients and conditions for regularity},
  booktitle = {Theory of Nonlinear Operators (Proc. Fourth Internat. Summer School, Berlin, 1975)},
  series = {Abh. Akad. Wiss. DDR, Abt. Math.-Natur.-Tech.},
  number = {1977, 1},
  publisher = {Akademie-Verlag}, address = {Berlin},
  pages = {197--206}, year = {1977}
}

@article{Oh93,
  author = {Oh, Yong-Geun},
  title = {Volume minimization of {Lagrangian} submanifolds under {Hamiltonian} deformations},
  journal = {Math. Z.},
  volume = {212}, pages = {175--192}, year = {1993},
  doi = {10.1007/BF02571651}
}

@article{SY00,
  author = {{\v{S}}ver{\'a}k, Vladim{\'i}r and Yan, Xiaodong},
  title = {A singular minimizer of a smooth strongly convex functional in three dimensions},
  journal = {Calc. Var. Partial Differential Equations},
  volume = {10}, number = {3}, pages = {213--221}, year = {2000},
  doi = {10.1007/s005260050151}
}

@article{SY02,
  author = {{\v{S}}ver{\'a}k, Vladim{\'i}r and Yan, Xiaodong},
  title = {Non-{Lipschitz} minimizers of smooth uniformly convex functionals},
  journal = {Proc. Natl. Acad. Sci. USA},
  volume = {99}, number = {24}, pages = {15269--15276}, year = {2002},
  doi = {10.1073/pnas.222494699}
}

@article{Mindlin64,
  author = {Mindlin, R. D.},
  title = {Micro-structure in linear elasticity},
  journal = {Arch. Rational Mech. Anal.},
  volume = {16}, pages = {51--78}, year = {1964},
  doi = {10.1007/BF00248490},
  url = {https://link.springer.com/article/10.1007/BF00248490}
}

@article{JK00,
  author = {Jin, Weimin and Kohn, Robert V.},
  title = {Singular perturbation and the energy of folds},
  journal = {J. Nonlinear Sci.},
  volume = {10}, number = {3}, pages = {355--390}, year = {2000},
  doi = {10.1007/s003329910014},
  url = {https://link.springer.com/article/10.1007/s003329910014}
}

@book{Morrey66,
  author = {Morrey, Jr., Charles B.},
  title = {Multiple integrals in the calculus of variations},
  series = {Die Grundlehren der mathematischen Wissenschaften},
  volume = {130},
  publisher = {Springer-Verlag}, address = {Berlin--Heidelberg--New York},
  year = {1966},
  doi = {10.1007/978-3-540-69952-1},
  url = {https://link.springer.com/book/10.1007/978-3-540-69952-1}
}

@article{GiMi68,
  author = {Giusti, Enrico and Miranda, Mario},
  title = {Sulla regolarit\`a delle soluzioni deboli di una classe di sistemi ellittici quasi-lineari},
  journal = {Arch. Rational Mech. Anal.},
  volume = {31}, number = {3}, pages = {173--184}, year = {1968},
  doi = {10.1007/BF00282679},
  url = {https://link.springer.com/article/10.1007/BF00282679}
}

@article{Giusti1969Singolarita,
  author   = {Giusti, Enrico},
  title    = {Precisazione delle funzioni di {$H^{1,p}$} e singolarit{\`a}
              delle soluzioni deboli di sistemi ellittici non lineari},
  journal  = {Boll. Un. Mat. Ital. (4)},
  volume   = {2},
  year     = {1969},
  pages    = {71--76},
  mrnumber = {0243190}
}

@article{Oh90,
  author = {Oh, Yong-Geun},
  title = {Second variation and stabilities of minimal {Lagrangian} submanifolds in {K\"ahler} manifolds},
  journal = {Invent. Math.},
  volume = {101}, pages = {501--519}, year = {1990},
  doi = {10.1007/BF01231513},
  url = {https://link.springer.com/article/10.1007/BF01231513}
}

@article{SW01,
  author = {Schoen, Richard and Wolfson, Jon G.},
  title = {Minimizing area among {Lagrangian} surfaces: the mapping problem},
  journal = {J. Differential Geom.},
  volume = {58}, number = {1}, pages = {1--86}, year = {2001},
  doi = {10.4310/jdg/1090348282},
  eprint = {math/0008244}, archivePrefix = {arXiv},
  url = {https://arxiv.org/abs/math/0008244}
}

@article{TW05,
  author  = {Trudinger, Neil S. and Wang, Xu-Jia},
  title   = {The affine {Plateau} problem},
  journal = {Journal of the American Mathematical Society},
  year    = {2005},
  volume  = {18},
  number  = {2},
  pages   = {253--289},
  doi     = {10.1090/S0894-0347-05-00475-3}
}

@article{Don02,
  author  = {Donaldson, S. K.},
  title   = {Scalar curvature and stability of toric varieties},
  journal = {Journal of Differential Geometry},
  year    = {2002},
  volume  = {62},
  number  = {2},
  pages   = {289--349},
  doi     = {10.4310/jdg/1090950195}
}

@article{Tou62,
  author  = {Toupin, R. A.},
  title   = {Elastic materials with couple-stresses},
  journal = {Archive for Rational Mechanics and Analysis},
  year    = {1962},
  volume  = {11},
  number  = {1},
  pages   = {385--414},
  doi     = {10.1007/BF00253945}
}

@article{Nec67,
  author  = {Ne{\v{c}}as, Jind{\v{r}}ich},
  title   = {Sur la r{\'e}gularit{\'e} des solutions variationnelles des
             {\'e}quations elliptiques non-lin{\'e}aires d'ordre {$2k$}
             en deux dimensions},
  journal = {Annali della Scuola Normale Superiore di Pisa,
             Classe di Scienze (3)},
  year    = {1967},
  volume  = {21},
  number  = {3},
  pages   = {427--457},
  url     = {https://www.numdam.org/item/ASNSP_1967_3_21_3_427_0/}
}

@article{FJM06,
  author  = {Friesecke, Gero and James, Richard D. and M{\"u}ller, Stefan},
  title   = {A hierarchy of plate models derived from nonlinear
             elasticity by {Gamma}-convergence},
  journal = {Archive for Rational Mechanics and Analysis},
  year    = {2006},
  volume  = {180},
  number  = {2},
  pages   = {183--236},
  doi     = {10.1007/s00205-005-0400-7}
}

@article{P,
    AUTHOR = {Paneitz, Stephen M.},
     TITLE = {A quartic conformally covariant differential operator for
              arbitrary pseudo-{R}iemannian manifolds (summary)},
   JOURNAL = {SIGMA Symmetry Integrability Geom. Methods Appl.},
  FJOURNAL = {SIGMA. Symmetry, Integrability and Geometry. Methods and
              Applications},
    VOLUME = {4},
      YEAR = {2008},
     PAGES = {Paper 036, 3},
   MRCLASS = {58J60 (58J70)},
  MRNUMBER = {2393291},
MRREVIEWER = {John Urbas},
       DOI = {10.3842/SIGMA.2008.036},
       URL = {https://doi.org/10.3842/SIGMA.2008.036},
}

@inproceedings{AG87,
  author    = {Aviles, Patricio and Giga, Yoshikazu},
  title     = {A mathematical problem related to the physical theory
               of liquid crystal configurations},
  booktitle = {Miniconference on Geometry/Partial Differential Equations, 2},
  editor    = {Hutchinson, John E. and Simon, Leon M.},
  series    = {Proceedings of the Centre for Mathematical Analysis},
  volume    = {12},
  publisher = {Australian National University},
  address   = {Canberra},
  year      = {1987},
  pages     = {1--16},
  url       = {https://maths.anu.edu.au/files/CMAProcVol12-AvilesGiga_2.pdf}
}

@article{Wil65,
  author  = {Willmore, T. J.},
  title   = {Note on embedded surfaces},
  journal = {An. \c{S}ti. Univ. ``Al. I. Cuza'' Ia\c{s}i
             Sec\c{t}. I a Mat. (N.S.)},
  volume  = {11B},
  year    = {1965},
  pages   = {493--496},
  url     = {https://www.math.uaic.ro/~annalsmath/pdf-uri_anale/remarkable-papers/Thomas-J.-Willmore-1965.pdf}
}

@article{GV15,
  author  = {Gursky, Matthew J. and Viaclovsky, Jeff A.},
  title   = {Rigidity and stability of {Einstein} metrics for
             quadratic curvature functionals},
  journal = {Journal f{\"u}r die reine und angewandte Mathematik},
  year    = {2015},
  volume  = {700},
  pages   = {37--91},
  doi     = {10.1515/crelle-2013-0024}
}

@article{DeFilippisKristensenMingioneStrong2026,
  author = {De Filippis, Cristiana and Kristensen, Jan
            and Mingione, Giuseppe and Strong, Aidan},
  title = {The singular set of minima of multiple integrals},
  year = {2026},
  note = {Preprint, arXiv:2609.11546},
  eprint = {2609.11546},
  archivePrefix = {arXiv},
  primaryClass = {math.AP},
  url = {https://arxiv.org/abs/2609.11546}
}

@article{GG78,
  author  = {Giaquinta, Mariano and Giusti, Enrico},
  title   = {Nonlinear elliptic systems with quadratic growth},
  journal = {Manuscripta Mathematica},
  volume  = {24},
  pages   = {323--349},
  year    = {1978},
  doi     = {10.1007/BF01167835}
}

@article{BW19,
  author  = {Bhattacharya, Arunima and Warren, Micah},
  title   = {Interior {Schauder} estimates for the fourth order
             {Hamiltonian} stationary equation in two dimensions},
  journal = {Proceedings of the American Mathematical Society},
  volume  = {147},
  number  = {8},
  pages   = {3471--3477},
  year    = {2019},
  doi     = {10.1090/proc/14484}
}

@article{Dimler2023,
  author  = {Dimler, Bryan},
  title   = {Partial regularity for Lipschitz solutions to the minimal surface system},
  journal = {Calc. Var. Partial Differential Equations},
  volume  = {62},
  year    = {2023},
  number  = {9},
  pages   = {260},
  doi     = {10.1007/s00526-023-02591-7},
}

@article{Duggan1986,
  author  = {Duggan, J. P.},
  title   = {Regularity theorems for varifolds with mean curvature},
  journal = {Indiana University Mathematics Journal},
  volume  = {35},
  number  = {1},
  pages   = {117--144},
  year    = {1986},
  doi     = {10.1512/iumj.1986.35.35006}
}

@misc{FiresterTsiamis26,
  author = {Firester, Benjy and Tsiamis, Raphael},
  title = {Topology of minimal surfaces in the sphere
           from capillarity},
  year = {2026},
  note = {Preprint, arXiv:2604.04928},
  eprint = {2604.04928},
  archivePrefix = {arXiv},
  primaryClass = {math.DG},
  url = {https://arxiv.org/abs/2604.04928}
}

@article{Wan16,
  author  = {Wang, Lizhou},
  title   = {Singular critical points for variational free boundary
             problem from isoparametric hypersurfaces},
  journal = {Communications in Contemporary Mathematics},
  volume  = {18},
  number  = {3},
  pages   = {1650010},
  year    = {2016},
  doi     = {10.1142/S0219199716500103}
}

@misc{FTZ26,
  author = {Firester, Benjy and Tsiamis, Raphael
            and Zhao, Zihui},
  title = {Singularity models for the {B}ernoulli free
           boundary problem from isoparametric hypersurfaces},
  year = {2026},
  note = {Preprint, arXiv:2609.31617},
  eprint = {2609.31617},
  archivePrefix = {arXiv},
  primaryClass = {math.AP},
  url = {https://arxiv.org/abs/2609.31617}
}

@article{CV2026,
  author = {Carducci, Matteo and Velichkov, Bozhidar},
  title = {A minimizing cone from {C}artan polynomials
           in the {A}lt--{C}affarelli problem},
  year = {2026},
  note = {Preprint, arXiv:2609.35179},
  eprint = {2609.35179},
  archivePrefix = {arXiv},
  primaryClass = {math.AP},
  url = {https://arxiv.org/abs/2609.35179}
}

@article{HW15,
  author  = {Hong, Guanghao and Wang, Lizhou},
  title   = {A singular critical point for some variational free
             boundary problem in five dimensions},
  journal = {Journal of Mathematical Analysis and Applications},
  volume  = {426},
  number  = {2},
  pages   = {855--863},
  year    = {2015},
  doi     = {10.1016/j.jmaa.2015.02.013}
}

@book{Weyl39,
  author    = {Weyl, Hermann},
  title     = {The Classical Groups: Their Invariants and Representations},
  publisher = {Princeton University Press},
  address   = {Princeton, NJ},
  year      = {1939}
}

@article{Moo20,
  author  = {Mooney, Connor},
  title   = {Minimizers of convex functionals with small degeneracy set},
  journal = {Calc. Var. Partial Differential Equations},
  year    = {2020},
  volume  = {59},
  number  = {2},
  pages   = {Paper No. 74, 19 pp.},
  doi     = {10.1007/s00526-020-1723-9}
}

@article{BB18,
  author  = {Bousquet, Pierre and Brasco, Lorenzo},
  title   = {{$C^1$} regularity of orthotropic {$p$}-harmonic functions in the plane},
  journal = {Anal. PDE},
  year    = {2018},
  volume  = {11},
  number  = {4},
  pages   = {813--854},
  doi     = {10.2140/apde.2018.11.813}
}
\end{document}